\title{ \bf Hyperbolic Structure of the\\ Equilateral Pentagon\\
{\large -- a quest of finding coordinates --}\\
}
\author{ \it J\"{u}rgen Richter-Gebert}

\documentclass[10pt]{article} 

\usepackage[letterpaper,
            bindingoffset=0.2in,
            left=.4in,
            right=1in,
            top=1in,
            bottom=1in,
            footskip=.25in]{geometry}
\usepackage{blindtext}
\usepackage{multicol}
\usepackage{lipsum}
\usepackage{cuted}
\usepackage{mwe}
\usepackage{graphicx}
\usepackage{wrapfig}
\usepackage[autostyle=true]{csquotes}

\usepackage{amsthm}
\newtheorem*{theorem*}{Theorem}
\newtheorem*{theoremA*}{Theorem A}
\newtheorem*{theoremB*}{Theorem B}
\newtheorem*{theoremC*}{Theorem C}

\newtheorem{theorem}{Theorem}

\newtheorem{definition}{Definition}

\usepackage{amssymb}
\usepackage[font=small,labelfont=bf]{caption} 
\usepackage{color}
\usepackage{mathtools}
\usepackage{amsmath}
\usepackage[dvipsnames]{xcolor}
\usepackage[all]{nowidow}
\usepackage{float}

\usepackage{soul} 
\usepackage[normalem]{ulem} 

\usepackage[color links=true, backref=page]{hyperref} 

\hypersetup{
  colorlinks,
  citecolor=blue,
  linkcolor=Red,
  urlcolor=Blue
  }

\begin{document}
\newpage
\maketitle





\def\upto{, \ldots ,}
\def\lines{\square}
\def\points{\bigcirc}
\def\lines{\vee}
\def\points{\wedge}
\parskip=1mm
\def\meet{\wedge}
\def\join{\vee}

\def\pts[#1]{\vee_{#1}}
\def\lns[#1]{\wedge_{#1}}
\def\coa#1{\color[rgb]{.2,.2,1}#1}
\def\cob#1{\color[rgb]{.7,0,.4}#1}
\def\coc#1{\color[rgb]{.8,0,0}#1}
\def\cod#1{\color[rgb]{.8,.6,0}#1}
\def\coe#1{\color[rgb]{0,.6,0}#1}

\setlength{\parindent}{.5cm}

%
%


\parskip=2mm
\begin{multicols}{2}

\section{Spaces of linkages}
It is one of the finest mathematical
moments to experience a deep relation between two seemingly unrelated 
objects. In particular, if one of the objects is of almost childish simplicity and the other one requires a certain mathematical depth and understanding to be comprehended. So it happened to me when a student of mine (Milan Arounas, whom I want to thank at this place for bringing that topic to the table) in a conversation mentioned that he recently heard about a relation of a pentagonal linkage and the regular hyperbolic tiling of type 
$(5,4)$. 
More precisely:  {\it The different combinatorial types of realisations of the equilateral pentagon can be associated with the vertices, edges and faces of the hyperbolic $(5,4)$-tiling} (compare \cite{AY98, CR07, KK19}). If in addition certain identifications  resembling the symmetry of the situation are made the correspondence is one-to-one and the globally underlying space (often called the moduli space of the linkage) takes the form of a orientable manifold of genus~4.


One of my main research fields is creating foundations and software for the visualisation of mathematical objects. So the instant, reflex like, response to Milan's statement was ``I want to write an interactive visualisation of that phenomenon''.  Writing interactive visualisations, in particular, if they should be continuous and smooth, often carries the challenge of finding concrete coordinates that create a position of a certain object dependent on a few control parameters. Very often the problems that arise in such a context are of mathematical nature rather than being computer science problems, since they require a meaningful connection from the control parameters to the derived object. Ideally, it is exactly this connection that provides mathematical insight to a potential user of the animation.

Let us formulate what one may consider the ultimate aim of the visualisation task in relation to the equilateral pentagon: {\it Provide an interactive animation where one moves one point in the hyperbolic plane  such that for every position a realisation of the equilateral pentagon is shown. At the corners and edges and face barycenters of a $(5,4)$-tiling the pentagon should agree with the corresponding combinatorial type that is associated with this position. Between those positions the animation should behave as smoothly as possible (ideally $C^\infty$ continuous, or even conformal in an appropriate sense).}

\begin{figure*}
\noindent
\begin{center}
\includegraphics[width=0.48\textwidth]{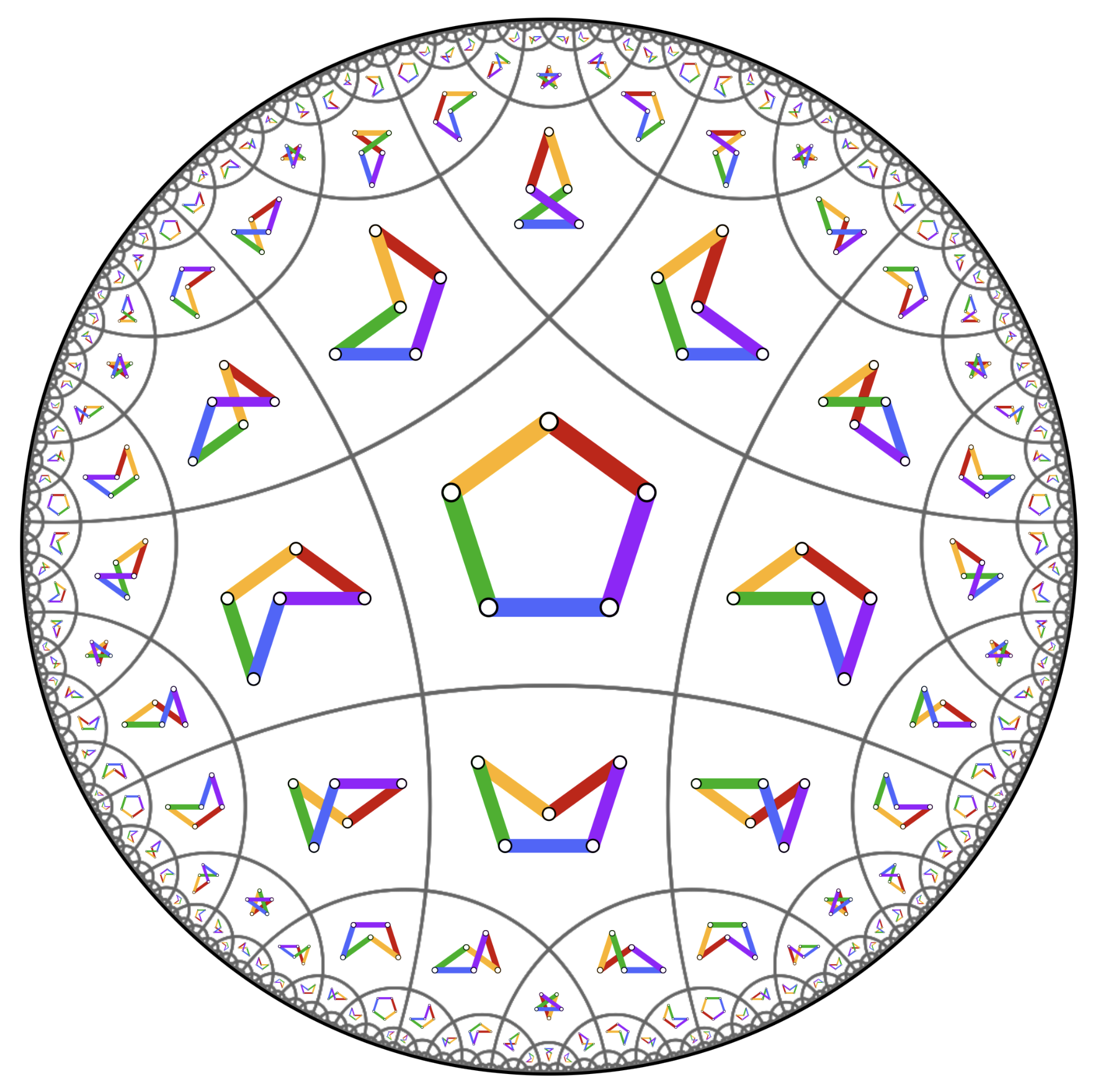}\quad
\includegraphics[width=0.48\textwidth]{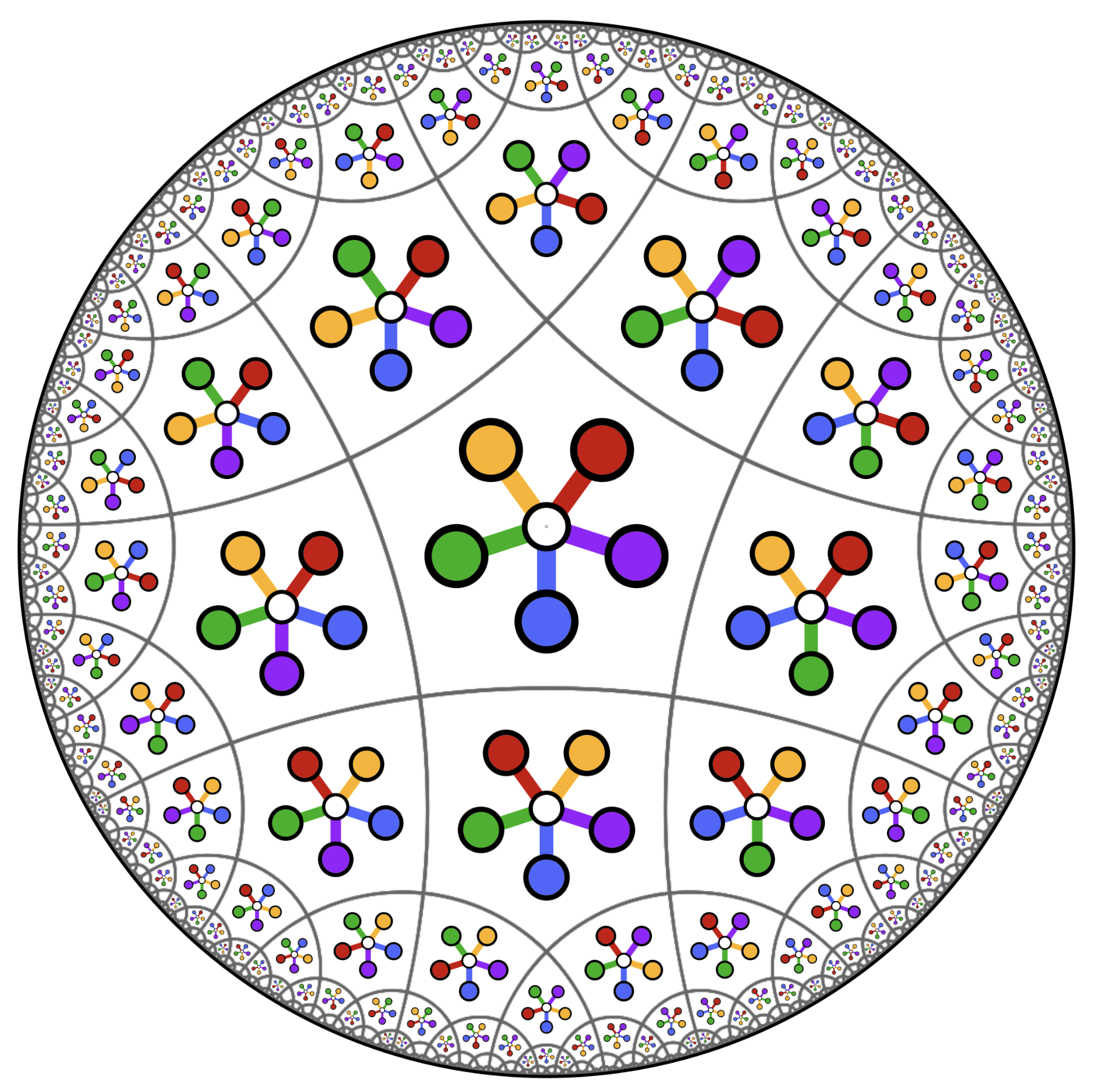}
\end{center}
\captionof{figure}{A first rough map of the pentagon space}~\label{fig:maps1}
\end{figure*}

\begin{figure*}
\noindent
\begin{center}
\includegraphics[width=0.95\textwidth]{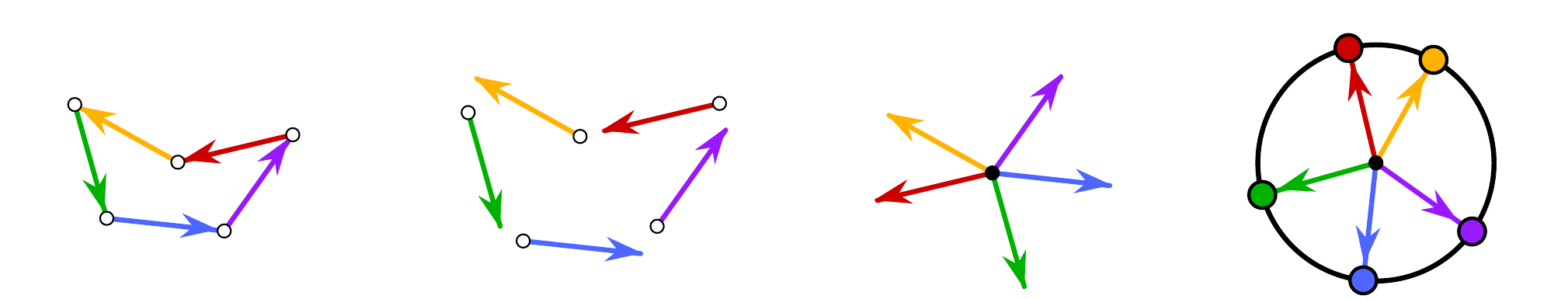}
\begin{picture}(0,0)
\put(-450,15){\footnotesize {$p_0$}}
\put(-410,12){\footnotesize {$p_1$}}
\put(-390,57){\footnotesize {$p_2$}}
\put(-426,49){\footnotesize {$p_3$}}
\put(-460,64){\footnotesize {$p_4$}}
\put(-310,7){\footnotesize {$v_0$}}
\put(-265,31){\footnotesize {$v_1$}}
\put(-282,60){\footnotesize {$v_2$}}
\put(-317,62){\footnotesize {$v_3$}}
\put(-342,36){\footnotesize {$v_4$}}
\put(-140,30){\footnotesize {$v_0$}}
\put(-158,71){\footnotesize {$v_1$}}
\put(-223,27){\footnotesize {$v_2$}}
\put(-219,57){\footnotesize {$v_3$}}
\put(-170,-3){\footnotesize {$v_4$}}
\put(-65,-6){\footnotesize {$v_0^\perp$}}
\put(-25,15){\footnotesize {$v_1^\perp$}}
\put(-75,87){\footnotesize {$v_2^\perp$}}
\put(-42,80){\footnotesize {$v_3^\perp$}}
\put(-113,28){\footnotesize {$v_4^\perp$}}

\end{picture}
\end{center}
\captionof{figure}{From pentagon linkages to 5-juzus.
Detach the edge vectors $v_i$, reassemble them around the origin,
rotate clockwise by $90^\circ$ and consider them as points on the unit circle.}~\label{fig:juzu}
\end{figure*}
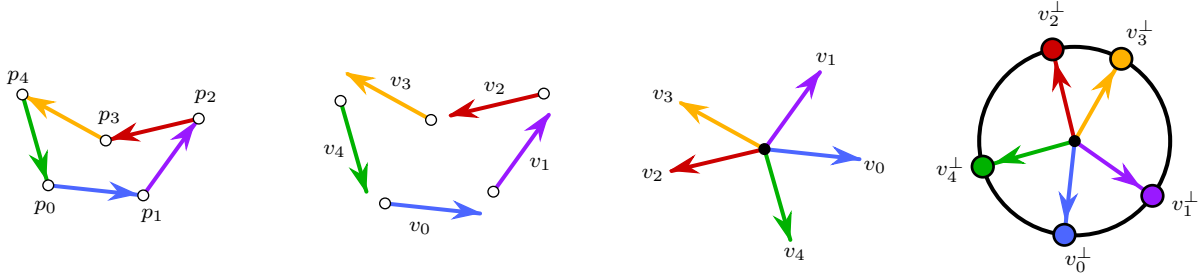

The left part of Figure~\ref{fig:maps1}  shows the corresponding equilateral pentagons for
the centers of the $(4,5)$-tiling. Each linkage shown there has uniform edge lengths. The edges are colored and along each pentagon the sequence of the colors is blue-purple-red-yellow-green. Lateron this sequence of colors will
correspond to the order of the edge vectors $(v_0,\ldots,v_4)$. Crossing the border of a pentagonal cell transforms a linkage into a neighboring one in which exactly two edge directions are interchanged. 
 On the right side of the image you see a related diagram.
There all edges are detached form the pentagon and  glued together at the (local) origin and rotated by $90^\circ$, such that the corresponding lines are orthogonal to the sides of the pentagons. Since in that representation all segments meet in the origin in further images we will  omit the edges themselves and only show the positions of the endpoints pointing away from the origin. These endpoints form a kind of bracelet consisting of five points on the unit circle in a certain cyclic order. They correspond to what  Yoshida and Apéry call 5{\it -juzus} in~\cite{AY98}. A very appropriate name:
 Juzus are Japanese prayer beads attached to a cyclic thread. We will use this name as well to refer to a configuration of five labelled points on the unit circle. 
  Figure~\ref{fig:juzu} explains the transition from pentagon linkage to a 5-juzu. 
 Like the pentagon linkages, juzus can be considered on several different abstraction levels.
We will distinguish between  a {\it juzu} (an element in  $(S^1)^5$) which will be 5 concrete points on the unit circle,
a  {\it pre-juzu}  (an element in  $(S^1)^5/S^1$) which is a juzu modulo rotation, and a  {\it combinatorial juzu} which encodes
the cyclic oder of a juzu.

In the right diagram of Figure~\ref{fig:maps1} it is relatively easy to 
understand the combinatorics of the entire situation. The $(5,4)$-tilling consists of right-angled hyperbolic pentagons. Four of them meet at each vertex of the tiling. Like in a checker board the sides align to (hyperbolic) lines (they are circular arcs in the Poincaré disk picture). 
Each line is associated to a pair of vectors that change their relative position
when crossing the line.
Since there are 5 colors we will have $10={5\choose 2}$ different types of lines.
The challenge of the animation is to associate a pentagon linkage to {\it each} point it the hyperbolic plane in a mathematically consistent and sound way that exhibits exactly this switching behaviour when crossing the lines in a continuous way.

In addition the animation should be  ``fluent'' this means that a mouse interaction results in a new image within milliseconds. This in addition may lead to a serious computer science challenges that will mainly be neglected here. You can find animations at 
\href{https://science-to-touch.com/Visu/Pentagon/Pentagon_Minimal.html}{\tt math-visuals.org/Pentagon}.


The method we present here provides a range of possible 
parameterisations. One of them will stand out by being conceptually very simple 
and at the same time easy to compute. 
A comparable but slightly different approach was independently developed by Lyle Ramshaw \cite{Ra22}.
There the inverse problem of finding hyperbolic coordinates for given pentagonal linkages was treated.
The parameterisation proposed there is covered in the range of parameterisation we present in Section 5.
While the approach in  \cite{Ra22} is more rooted in differential geometry the treatment in this article is based on projective considerations.



\section{The equilateral Pentagon}
Rather than working in the real plane $\mathbb{R}^2$ we will consider our points $p_i$ in the complex plane $\mathbb{C}$. 
 In what follows, indices will always be counted modulo~$5$, if not otherwise stated.

\begin{definition}
We call five points  $p_0,p_1,p_2,p_3,p_4\in \mathbb{C}$ that satisfy 
\[|p_0,p_1|=|p_1,p_2|=|p_2,p_3|=|p_3,p_4|=|p_4,p_0|=1\]
a realisation of the equilateral pentagon.
\end{definition}

\begin{figure}[H]
\noindent
\begin{center}
\includegraphics[width=0.45\textwidth]{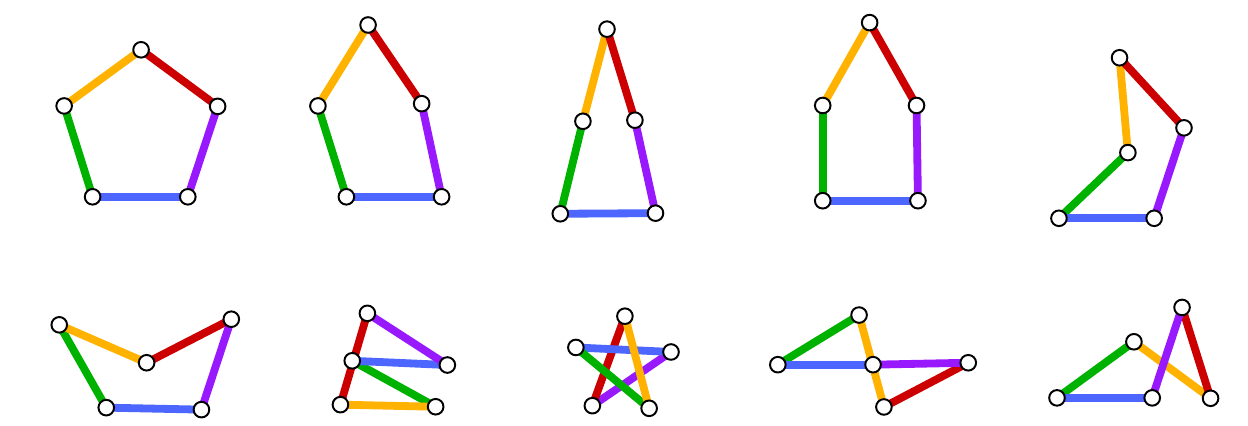}
\end{center}
\captionof{figure}{Some examples of equilateral pentagons}~\label{fig:examples}
\end{figure}

 Obviously the regular pentagon with side length~$1$ is such a realisation,
but there are many others.  Figure~\ref{fig:examples} shows a few examples. Each of these realisations can be considered do be assembled from five unit length segments forming a cycle. Due to the flexibility of the arrangement we will call it a equilateral pentagon {\it linkage}. In what follows we will drop the word equilateral, since all pentagon linkages we  consider will have uniform side length.
In an ideal world, where the bars can freely move through each other without mechanical hinderance we may think of the set of all realisations being realised by a hinge and bar mechanism of five unit length bars that cyclically close.

We  will mainly  consider differences of consecutive points
$v_k=p_{k+1}-p_{k}$.
Despite the fact that $v_k$  is also a complex number we will emphasise its role by calling it a {\it vector}.
Considering the $v_i$ will automatically mod-out translations. Pentagons that only differ by translation become the same object.
The choice of the five vectors $v_k$ is not completely arbitrary they have to satisfy 
\[
\sum v_k=0.
\]
For aesthetic  and perceptional reasons related to our images we make a further twist and rotate these vectors by $90^\circ$
counterclockwise and very often consider $v_k^\perp=-i\cdot v_k$ instead. Those are the points we draw in the juzus.

Taking into account that each $v_k$ is a point on the unit circle there are 5 real degrees of freedom for the choice of $v_k$. Two of them, are used to fix the vector sum property. If we are only interested in the {\it shape} of the pentagon we may furthermore identify pentagons that arise from each other by rotation. We can do so by setting $v_0=1$.
\begin{definition}
A normalised equilateral pentagon linkage (pentagon linkage, for short) is given by unit complex numbers $(v_0,\ldots,v_4)\in (S^1)^5$ with $v_0=1$ and $\sum v_k=0$. 
\end{definition}



 Locally, the space of pentagon linkages is a two dimensional manifold. We will refer to that space  by the letter $\mathcal{P}$.
  From the five original degrees of freedom two have been used for the closing condition and one is used for the normalisation of $v_0$. We will call this space  $\mathcal{P}$ of shapes the {\it realisation space} of the pentagon linkages (at other places this is called the moduli space).  In what follows we are doing no more and no less than asking for good parameterisations of  $\mathcal{P}$ by a point moving in $\mathbb{H}$. 
%

\begin{figure}[H]
\noindent
\begin{center}
\ \includegraphics[width=0.2\textwidth]{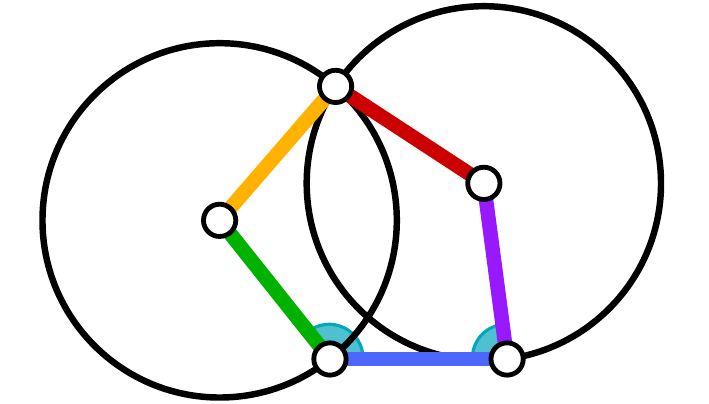}\quad
\includegraphics[width=0.22\textwidth]{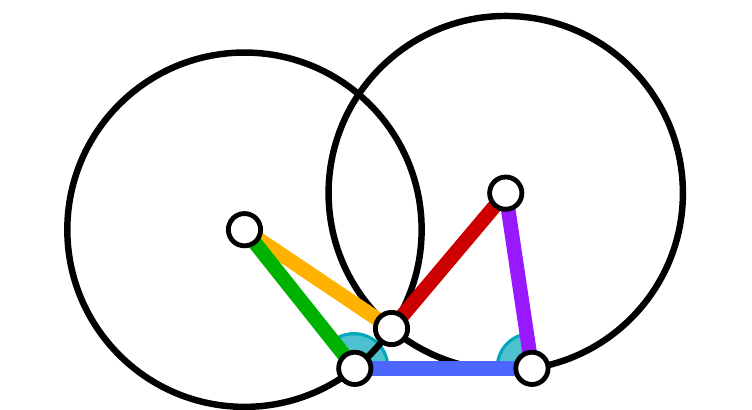}
\vskip2mm
\includegraphics[width=0.22\textwidth]{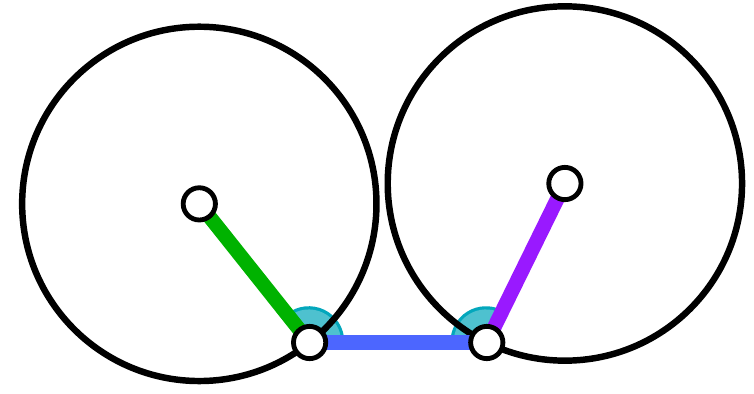}
\includegraphics[width=0.22\textwidth]{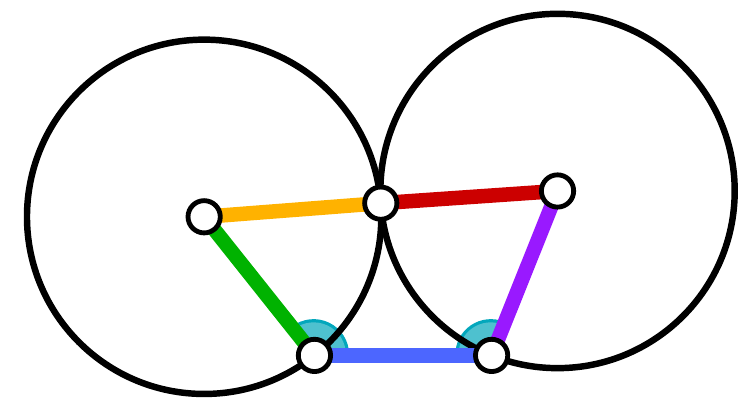}
\end{center}
\captionof{figure}{Naive parameterisation by two angels. The first two images  show the two solutions for the same parameters. The other two show a situation with no and a situation with a unique solution.}~\label{fig:naive}
\end{figure}

\noindent
\underline{\it A naive parameterisation:}
A first naive way to parameterise the realisation space is  close to its mechanical interpretation. We  fix $p_0=0, \ p_1=1$. Then we directly take
the angles $\alpha$ and $\beta$ that determine $v_1=e^{i\alpha}$ and
$v_4=e^{i\beta}$ they determine the positions of $p_2$ and $p_4$. The position of $p_3$ now must be such that it has unit distance to $p_2$ and $p_4$. This problem in general has two solutions that can be derived by intersecting two unit circles around $p_2$ and $p_4$. Transferred to real 2-dimensional coordinates this leads to a quadratic equation that has either  
two solutions (if $|p_2,p_4|<2$), one solution (if $|p_2,p_4|=2$) or has no solution at all if $p_2$ and $p_4$ are too far apart which may happen for certain choices of angles (see Figure~\ref{fig:naive}).

Such a straight forward coordinatisation has several disadvantages.  The obvious one is that for certain choices of parameters there are no solutions at all and for almost all  others there are two solutions. However, the perhaps biggest disadvantage is that it destroys the symmetry that is inherent to the problem. In this parameterisation the two points $p_0$ and $p_1$ play a role that is different from the other points. Also $p_3$ plays a special role. We will come back to the issue of choosing a better realisation soon. Before that we will deal with the combinatorial structure of the parameter space.

\section{The combinatorics}
We  now study how the realisation space $\mathcal{P}$
can be decomposed into reasonable cells  reflecting the combinatorial type of a realisation.

\subsection{A combinatorial manifold}
We have seen that five vectors  $(v_0,\ldots,v_4)\in (S^1)^5$ that satisfy
$\sum{v_k}=0$ determine the shape of a pentagon. 
  By the commutativity of vector addition every permutation $(\pi(v_0),\ldots,\pi(v_4))$ also satisfies the summation property and gives rise to some other realisation of the pentagon linkage with a different cyclic order of the points on the juzu.
We will now classify the different combinatorial types of such shapes by the combinatorics of the juzus.

If we orient the $S^1$ and require  $v_0,\ldots,v_4$ to be distinct there are 
 ${5!\over 5}={120\over 5}=24$ different possibilities for the cyclic order. They correspond to  {\it non-degenerate} situations. The non-degenerate pentagon linkages of the same type topologically form an open disk within~$\mathcal{P}$.
Our next observation is that in a proper pentagon it may not happen that three  vectors coincide. In that case the vector sum of these three vectors would have an absolute value of $3$ which cannot be compensated by the other two. 



However, it may happen that two vectors coincide. There are ${5 \choose 2}\cdot 3!=60$ combinatorial possibilities for that to happen. Any pair may coincide and the remaining points can be ordered in $3!$ possible ways. We will call such a situation {\it singly degenerate}.  The singly degenerate pentagon linkages of the same type topologically form an open segment within~$\mathcal{P}$, since one degree of freedom is needed to keep the two vectors aligned.
Every non-degenerate situation is surrounded by 5 possible singly degenerate situations since every pair of adjacent vectors may become coincident.

\begin{figure}[H]
\centering
\includegraphics[width=0.45\textwidth]{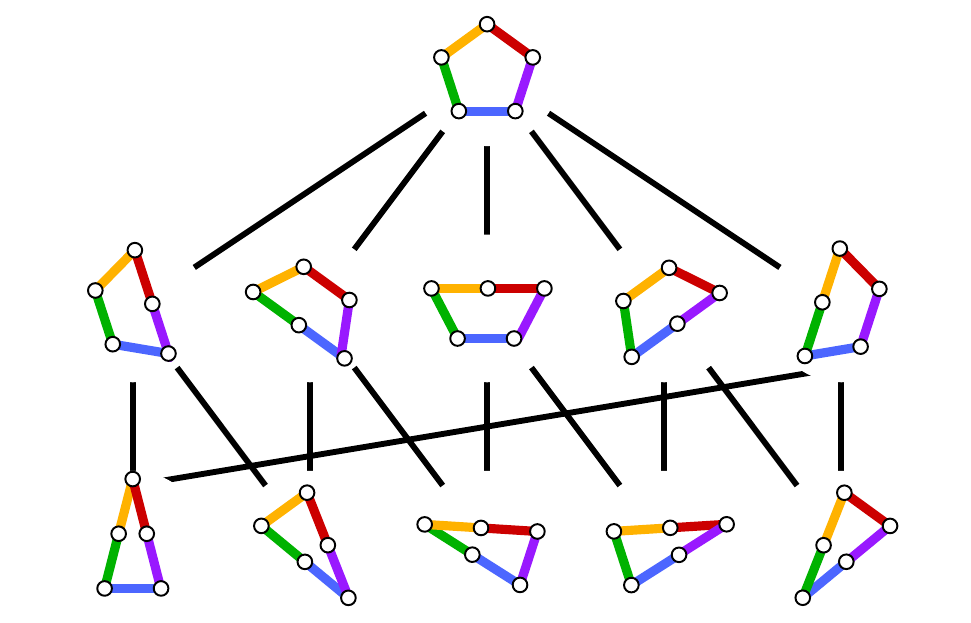}
\caption{\label{fig:partial}A natural partial order}
\end{figure}

\begin{figure*}[t]
\noindent
\begin{center}
\includegraphics[width=0.43\textwidth]{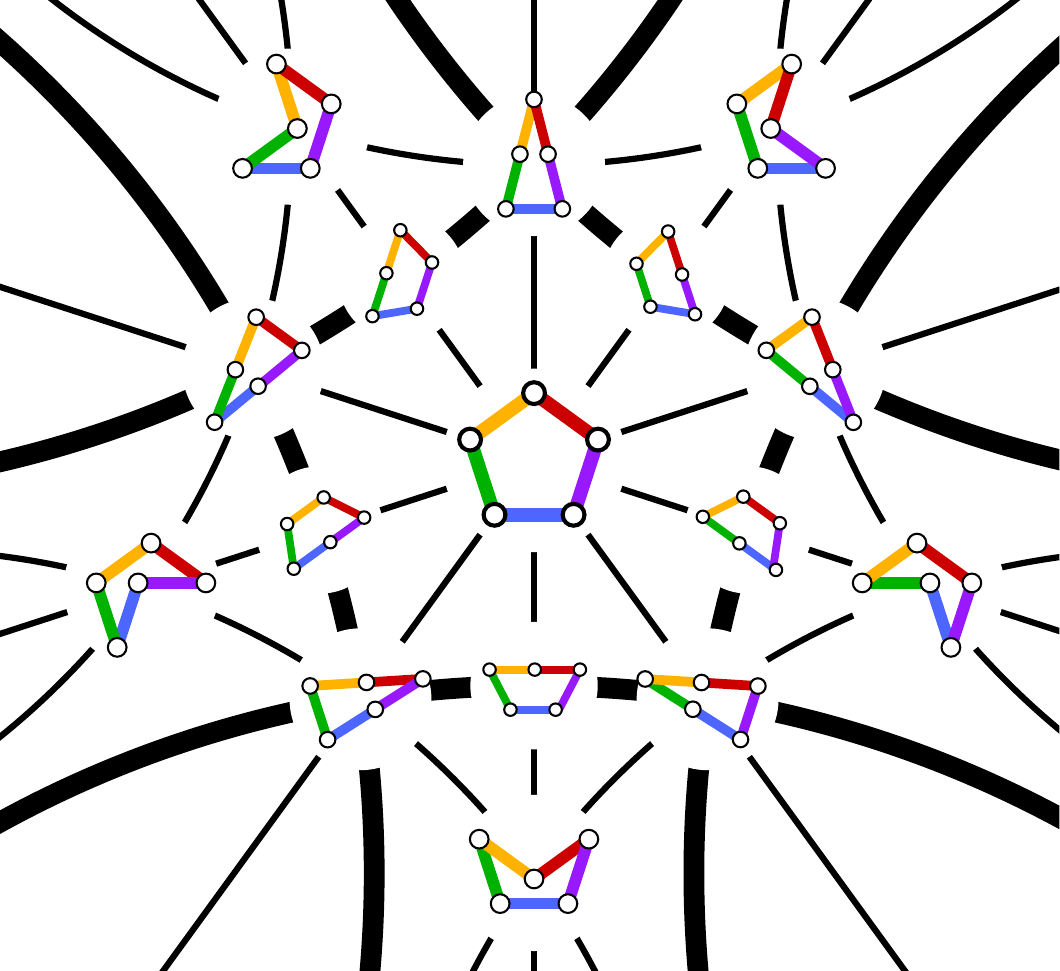}\qquad
\includegraphics[width=0.43\textwidth]{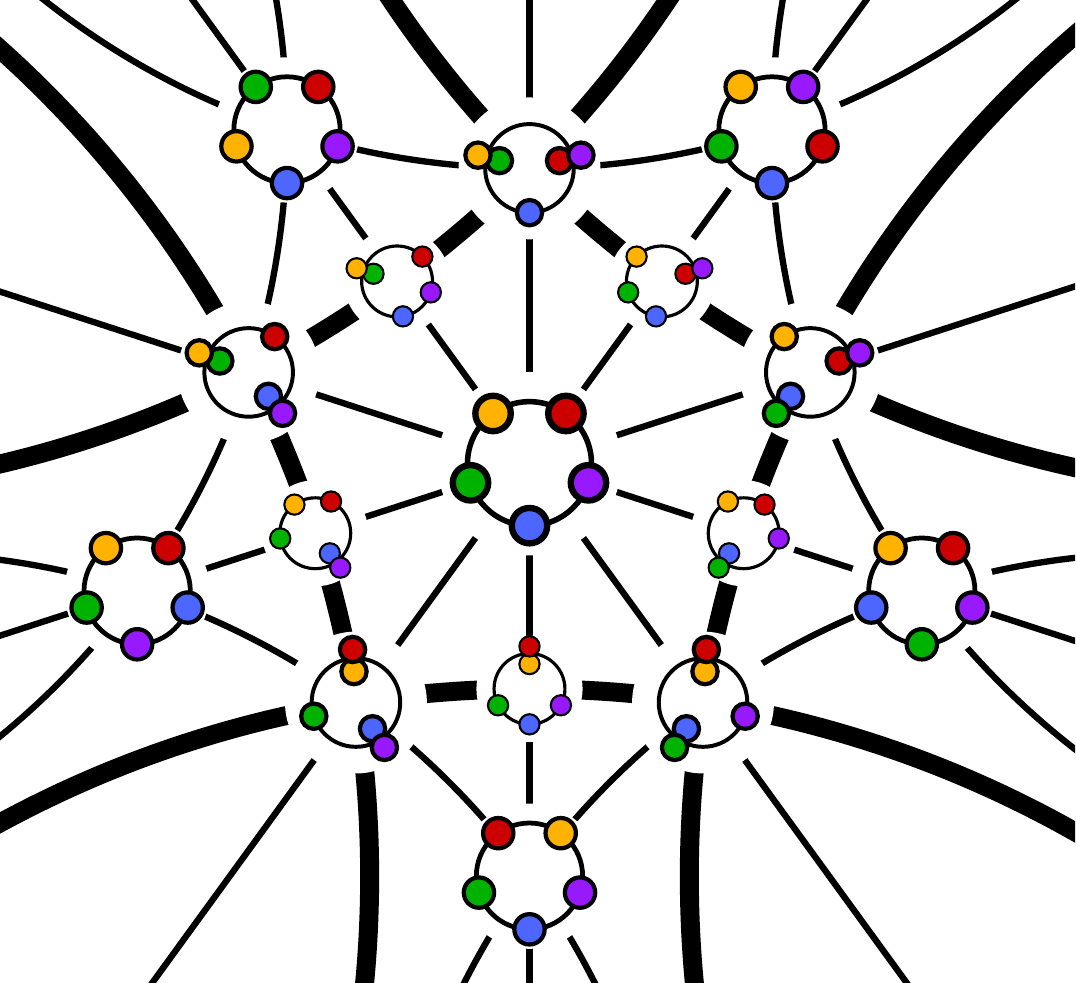}
\end{center}
\captionof{figure}{The neighbourhood of the regular pentagon}~\label{fig:maps2}
\end{figure*}

There are also {\it doubly degenerate} situations in which two independent pairs of vectors coincide. In the realisation space $\mathcal{P}$ they will be resembled by certain points. In terms of shapes of pentagons they may correspond to the situation in which two pairs of adjacent edges become collinear and the pentagon becomes a triangle with side lengths $2$, $2$ and $1$. However, there are  other situations where the aligned vectors are not adjacent in the sequence of edges. For an example of this see the second picture in the second row of Figure~\ref{fig:examples} which is double degenerate as well. A simple counting argument shows that there are ${5 \choose 2}\cdot {3 \choose 2} =30$ possibilities for double degeneracy to happen. Those $24+60+30$ different situations are the $114$ possible {\it combinatorial types} of the equilateral pentagon.

There is a natural graded partial order
by the amount of degeneration that is present in a combinatorial type.
We set  $t_1<t_2$ if $t_2$ can be reached from $t_1$ by moving coinciding vectors apart. Thus a singly degenerate type has two non-degenerate types that are smaller than it. 
The Hasse diagram in Figure~\ref{fig:partial}  exemplifies this partial order by 
showing all elements below our standard regular pentagon. The Hasse diagram itself  corresponds to the face lattice of a pentagon. Analogously the sub-lattice of all elements above a doubly degenerate situation correspond to a four-valent point surrounded by four faces.

\subsection{A hyperbolic tesselation}

Let us now turn to the global structure of $\mathcal{P}$.
Figure \ref{fig:maps2} shows the regular pentagon and a halo of its immediate neighborhood. It is worth contemplating a bit on this image. The geometric position of the shown pentagon linkages is chosen in a way that relates to their position in the corresponding tiling. The big black circular arcs in the background are hyperbolic lines consisting of the edges of the hyperbolic $(5,4)$-tiling. The interior points of those edges correspond to the singly degenerate situations.
Along each such line one pair of edges of the linkage points into the same direction. At places
where two such lines intersect (the vertices of the $(5,4)$-tiling)  {\it two} pairs of edges point in the same direction.
These are the doubly degenerate situations.

\begin{figure*}[t]
\noindent
\begin{center}
\includegraphics[width=0.67\textwidth]{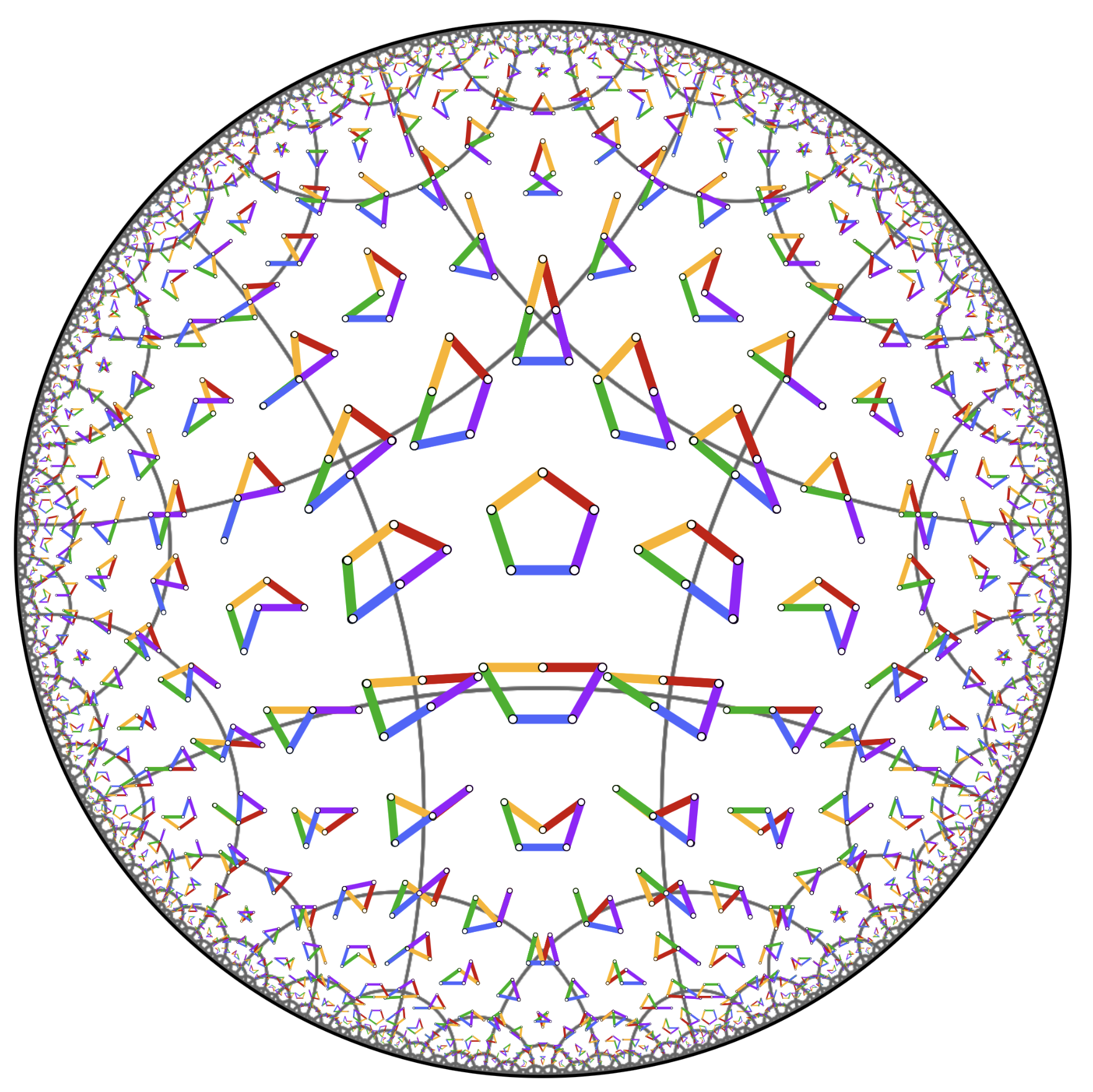}\qquad
\end{center}
\captionof{figure}{The pentagon linkages at the centres of the cells in each dimension.}~\label{fig:mapsFull}
\end{figure*}

The rotations of the pentagon linkages have been chosen in a very specific way.
This becomes  transparent if we again consider the 5-juzus realised by the $v_k^\perp$ that are shown in the right part of Figure \ref{fig:maps2}.
Consider the {\it thin} black symmetry lines in the background  of the image. 
Along those lines the juzus all support a mirror symmetry 
that have a fixed point of the same color. For instance, along the  vertical 
main axes of symmetry all juzu images are  vertically symmetric and have the blue point as fixed point.
It will turn out later that this can be done consistently  throughout the entire tiling.

\medskip

%
%
This situation is considerably well behaved.
The non-degenerate types correspond to pentagons. The doubly degenerate types correspond to vertices of degree four.
These local situations give rise to a two-dimensional structure in a cell complex and we can stitch all types together to obtain a combinatorial 2-manifold  describing the types of $\mathcal{P}$. 
Calculating the Euler characteristic $V-E+F$ we get
\[
\chi(\mathcal{P})=24-60+30=-6
\]
and shows that the manifold has genus 4.
Taking everything together we get

\begin{theorem} \label{thm:space}
The realisation space $\mathcal{P}$ has the structure of a 2-dimensional manifold of genus 4. This manifold is dissected into a  cell complex consisting of 24 pentagons such that four of them meet at each vertex.
\end{theorem}

The regular $(5,4)$-tiling as seen  Figure \ref{fig:maps1} perfectly fits the local situation on $\mathcal{P}$.
By suitably cutting out a 
finite region of 24 pentagons of the $(5,4)$-tiling and identifying corresponding parts of the boundary we  obtain a  representation of $\mathcal{P}$. 
%
%

Alternatively (and this is the way we will prefer in this article) we may leave the  $(5,4)$-tiling unchanged by covering it in a regular way with copies of the compound of 24 pentagons. This is analogous to representing a torus by an euclidean  pattern that repeats in two directions.

By that every point in the hyperbolic plane gets associated with exactly one element of $\mathcal{P}$. Thus we can navigate the space of equilateral pentagons by moving a point in the hyperbolic plane $\mathbb{H}$. Finding such a perfectly symmetric parametrisation that is as smooth as possible is the ultimate goal of the visualisation challenge. 
Figure \ref{fig:mapsFull} shows all linkages that we have nailed down to certain positions by our considerations yet. In what follows we will associate the center of that central tile (the origin) with the most natural representation of the linkage: The regular pentagon defined by the edge vectors $v_k=e^{k\cdot {2i\pi\over 5}}$. 


\section{\!A democratic parameter space}
Now  we go beyond combinatorial types 
to concrete parameterisations. 
We are looking at symmetric parameter space which reflects the symmetry of the situation.

In his brilliant book ``Hyperbolic Functions, my love'' Masaaki Yoshida describes the concept of a ``democratic parameterisation'' of a symmetric object  \cite{Y97}. Roughly speaking, in a democratic parameterisation the inherent symmetry of a class of objects should be reflected by the symmetry of the parameterisation. As he says in his book ``all points should be treated equally''. In our situation we want to treat each of the vectors $v_0,\ldots v_4$ in exactly the same way. 

While the space $(S^1)^5$ is perfectly symmetric in our five vectors, we have the problem that not every configuration of 5 vectors in $(S^1)^5$ represents a proper pentagonal linkage.  The five vectors  have to satisfy the summation property $\sum v_i=0$. 
However, there is a nice and elegant way to bypass this problem. It is based on a theorem by Boris Springborn \cite{Sp05} from wich we need a special case. The theorem is about normalisation of a point configuration via Möbius transformations, such that the center of gravity is moved to the origin. In its original form Springborns Theorem refers to $n$ {\it distinct} points on the unit sphere. But the same method also applies, if there are at least $d+2$ distinct points on a $d$-sphere (private communication). An independent version of this statement that also refers to continuous mass distributions  can be found in~\cite{KM95}, Theorem~4.
We only state the special case relevant to our application (Springborn's original result states that a similar theorem holds in every dimension and for every number of points)

%
%
\begin{theorem}
Let $(v_0,\ldots v_4)\in(S^1)^5$ be an collection of points on the complex unit circle such that no three points coincide. Then there is a unique orientation preserving Möbius transformation $\mu$ that leaves the circle invariant such that the center of  of gravity ${1\over 5}\sum\mu(v_k)$ equals the origin and $\mu(v_0)=1$\end{theorem}
%


Springborn's Theorem in that form is almost tailor-made for our purposes.
If we choose arbitrary five vectors $V=(v_0,\ldots v_4)$ on the unit circle the theorem associates a {\bf unique} shape of a pentagon linkage to $V$. Moreover, if a set $W=(w_0,\ldots w_4)$ of points on the unit circle differs from $V$ by a Möbius transformation it results in the same shape of a Pentagon linkage. In other words we parameterised the space of Pentagon linkages by 
five points on the unit circle none of which plays a special role. {\it  This is very democratic!}

It is a good exercise to see how counting degrees of freedom reflects the situation. Via stereographic projection the unit circle is isomorphic to the real projective line $\mathbb{RP}^1$. Möbius transformations on $S^1$ correspond to projective transformations on $\mathbb{RP}^1$. Such projective transformations have three real degrees of freedom. Thus up to Möbius transformation the five points are controlled by $5-3=$ degrees of freedom, as it should be the case for our pentagon linkage.

By all these  consideration we silently  have identified the space $\mathcal{P}$
with two other interesting objects. The space of five points on $S^1$ no three of which are identical, modulo Möbius transformations. And the space of
five points in $\mathbb{RP}^1$ no three of which are identical, modulo projective transformations. The five points in both representations play an absolutely democratic role. None of them is special.
%
%
\section{Exploiting symmetry}
Now we enter the challenging part of the problem. Figure \ref{fig:mapsFull} gives a rough understanding how a continuous parameterisation should behave.
Roughly speaking, we are looking for a parameterisation in which the symmetries of the hyperbolic tiling  find their counterparts in symmetries of the realisation space $\mathcal{P}$. For that reason let us first understand how symmetries of the $(5,4)$-tiling interact with $\mathcal{P}$ for those parts we know from Figure \ref{fig:mapsFull}.
We will focus on the juzu representation of the linkages and consider very concrete realisations of the juzus. Although in the realisation space $\mathcal{P}$ we identify linkages that differ only by rotation, in this chapter we will {\it explicitly not do so}. 
We will consider five points on the unit circle without further identification: $(v_0\upto v_4)\in (S^1)^5$.

In a sense our strategy will be very bold.
We will directly use the symmetries that can be observed from the
juzu image in Figure \ref{fig:maps2} and formulate general requirements from them. Then we will show that if we apply Springborn's normalisation to a parameterisation that satisfies these requirements  we get a suitable solution of our visualisation challenge. Finally, we will single out one specific way that can be computed easily.

%
%

\subsection{Symmetries of the tiling}
First let us list the symmetries within the $(5,4)$-tiling. 
The symmetry of the tiling is a (hyberbolic) reflection group generated by a triangle with corner angles $90^\circ, 45^\circ, 36^\circ$. In our drawings these reflections are generated by two types of hyperbolic lines, lets call then $A$ and $B$. 
The axes of Type~$A$ reflections are the lines supported by the edges of the pentagons (the lines shown in Figures \ref{fig:maps1} and \ref{fig:mapsFull}). The axes of Type $B$ reflections are the symmetry lines of the pentagons (some of them you see in Figure \ref{fig:maps2} as thin black lines; both types of symmetry lines can also be seen in
Figure~\ref{fig:hypreflections}  and Figure~\ref{fig:tiling1}).

Those two types of reflections play  different roles with respect to our pentagon linkages and juzus.  Type $B$ reflections correspond to reflections of the vectors
$(v_0,\ldots,v_4)$ in an axes though the origin.
One might have the impression that  Figure \ref{fig:mapsFull} only supports one such reflection symmetry of Type~$B$ (the central bilateral axis). However, this only comes from the specific angle of rotation that we have chosen for each of the pentagon linkages. If we rotate all of them  by $k\cdot 72^\circ$ simultaneously  other central reflection axes are supported. 
For pentagon linkages that are associated to positions on  those reflection lines the vector configuration $(v_0,\ldots,v_4)$ 
exhibits a mirror symmetry.

The action of type $A$ reflections is somewhat simpler. They create transpositions of two vectors in 
$(v_0,\ldots,v_4)$. Thus along the corresponding reflection lines two vectors will coincide.  Both effects can be observed in Figure \ref{fig:maps2}.

\subsection{Towards coordinates}

We are now heading towards a concrete map
$\Psi:\mathbb{H}\to (S^1)^5$
that assigns to a point $p\in\mathbb{H}$ a concrete list of five
juzu vectors $(v_0^\perp,\ldots,v_4^\perp)$. 

\begin{figure}[H]
\noindent
\begin{center}
\includegraphics[width=0.4\textwidth]{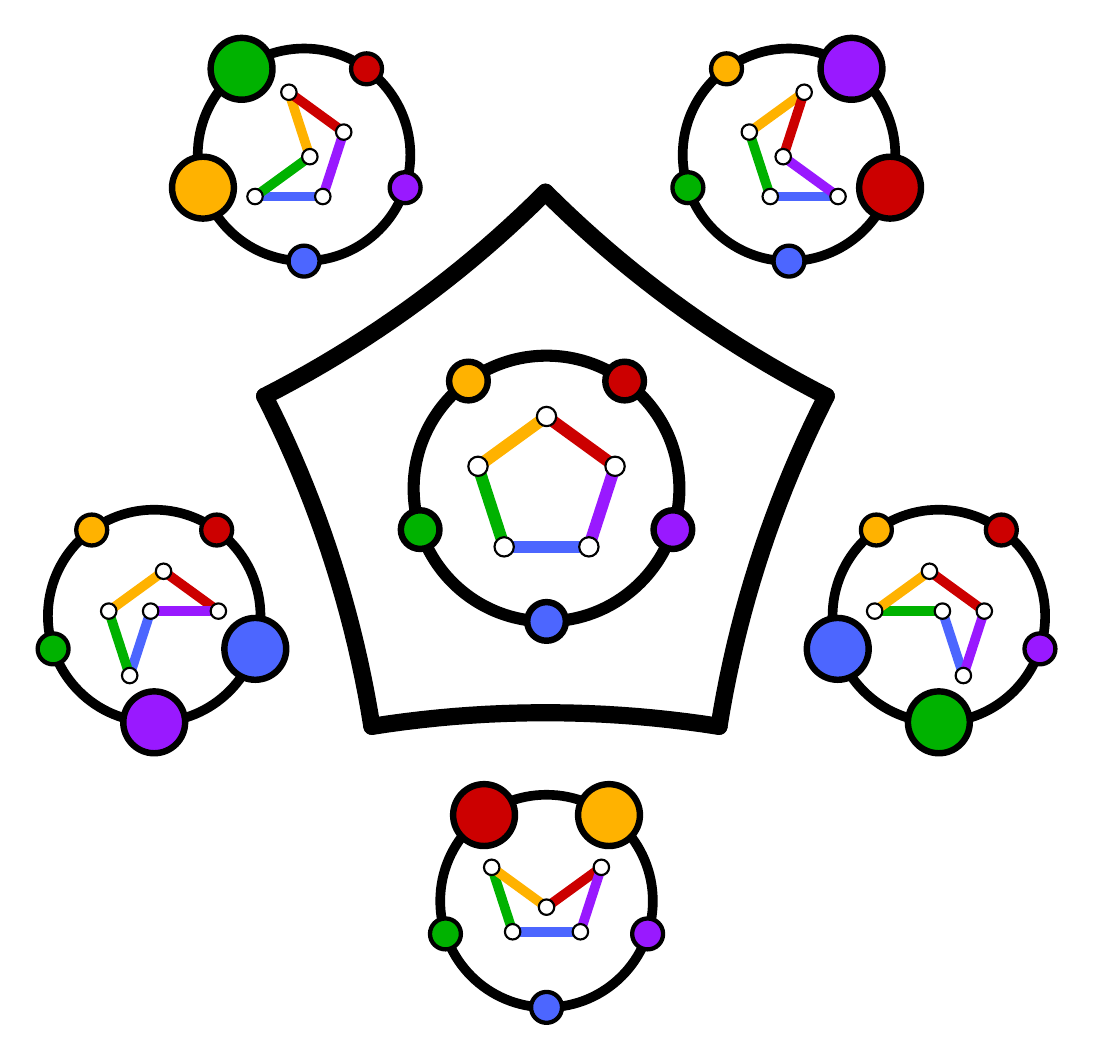}
\end{center}
\captionof{figure}{Neighbours of a pentagon}~\label{fig:rotate}
\end{figure}

There is one  subtlety that we have to keep in mind when comparing the symmetries in $\mathcal{P}$ and the symmetries of the $(5,4)$-tiling. Consider a hyperbolic  pentagon $X$ of the tiling and the juzu vectors $V=(v_0^\perp,\ldots, v_4^\perp)$ associated to the barycenter of $X$. We may assume that the order is chosen in a way that they correspond to a regular pentagon. Figure~\ref{fig:rotate} 
illustrates the situation.
When we cross an edge of $X$ 
adjacent vectors in $V$ will change places. The image highlights those  vectors that switched their roles by making them slightly bigger. 
If you look at the highlighted points of the  surrounding juzus they seem to travel around the
circle twice when making a full tour around $X$. This effect is the reason why we cannot produce a picture (neither of linkages nor of vectors or juzus) that reflects all mirror symmetries of the tiling simultaneously.

Let us now assume that we have 
a concrete function $\Psi:\mathbb{H}\to (S^1)^5$
that associates a  list of 5 juzu vectors to every single point in the hyperbolic disc. We will collect properties that such a function should satisfy. For that we assume that $\mathbb{H}$ is embedded into the complex plane via the Poincaré disc model bounded by the unit circle. From now on we set $\alpha:={2\pi/ 5}$.

\noindent
{{\bf Requirement 0} (Normalisation):} 
Let us start with a normalisation  condition:
For the origin $0\in\mathbb{H}$ we want 
\[
\Psi(0)=-i\cdot(e^{0},e^{ i\alpha },e^{2 i\alpha },e^{3 i\alpha },e^{4i\alpha } )\]
with that choice the origin corresponds to the regular pentagon. 
Then factor $-i$ resembles the rotation by $90^\circ$ to get the juzu vectors from the vectors.

\begin{figure}[H]
\noindent
\centering
\includegraphics[width=0.45\textwidth]{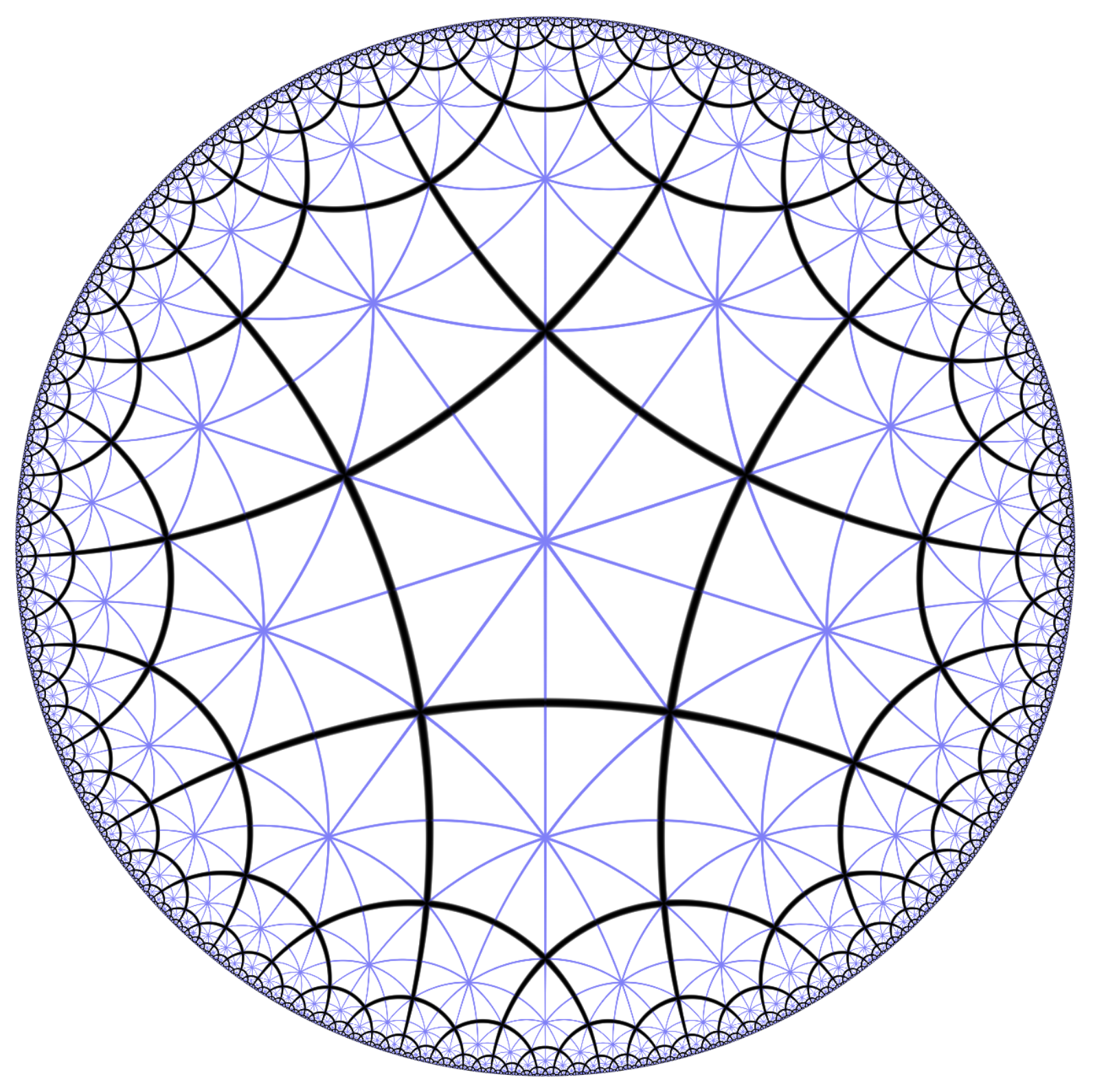}
\begin{picture}(0,0)
\put(-120,225){\footnotesize {$r_0$}}
\put(-186,206){\footnotesize {$r_1$}}
\put(-230,150){\footnotesize {$r_2$}}
\put(-230,78){\footnotesize {$r_3$}}
\put(-188,17){\footnotesize {$r_4$}}
\put(-22,50){\footnotesize {$s$}}
\end{picture}
\caption{\label{fig:tiling1}The central reflections}
\end{figure}

The full symmetry group of the $(4,5)$-tiling is
generated by three reflections in the sides of a triangle with hyperbolic angles $(90^\circ,45^\circ,36^\circ)$. 
Those triangles form the barycentric subdivision of the pentagons in the tiling. 
For each such triangle the mirror axes supported by the three sides generate the reflection group. Figure~\ref{fig:tiling1} shows how the $(5,4)$-tiling is regularly decomposed into 10 triangles per pentagon. 
Let $s,r_0,r_1$ be hyperbolic reflections corresponding to the sides of such a triangle. Each reflection is an involution $\mathbb{H}\to \mathbb{H}$ and together they generate this group. The underlying group can be described in a purely 
combinatorial way without explicitly referring to the geometric situation.
For this start with a free group generated by $s,r_0,r_1$ and consider the quotient group modulo the relations
\[
s^2=r_0^2=r_1^2=(sr_0)^2=(sr_1)^4=(r_0r_1)^5=1.
\]
The resulting group will be called $\mathcal{R}$.
Notice that it acts on $\mathbb{H}$.
 Along the boundary, Figure~\ref{fig:tiling1}  also indicates the position of the central reflections and the position of the pentagon side $s$. We also label the other central reflections $r_3,r_4,r_5$.
In the coordinate system we have chosen $r_k$ may be expressed as
\[
r_k(z):=-e^{i\cdot k\alpha}\cdot {\overline{z}}.
\]

\begin{figure*}[t]
\noindent
\begin{center}
\includegraphics[width=0.60\textwidth]{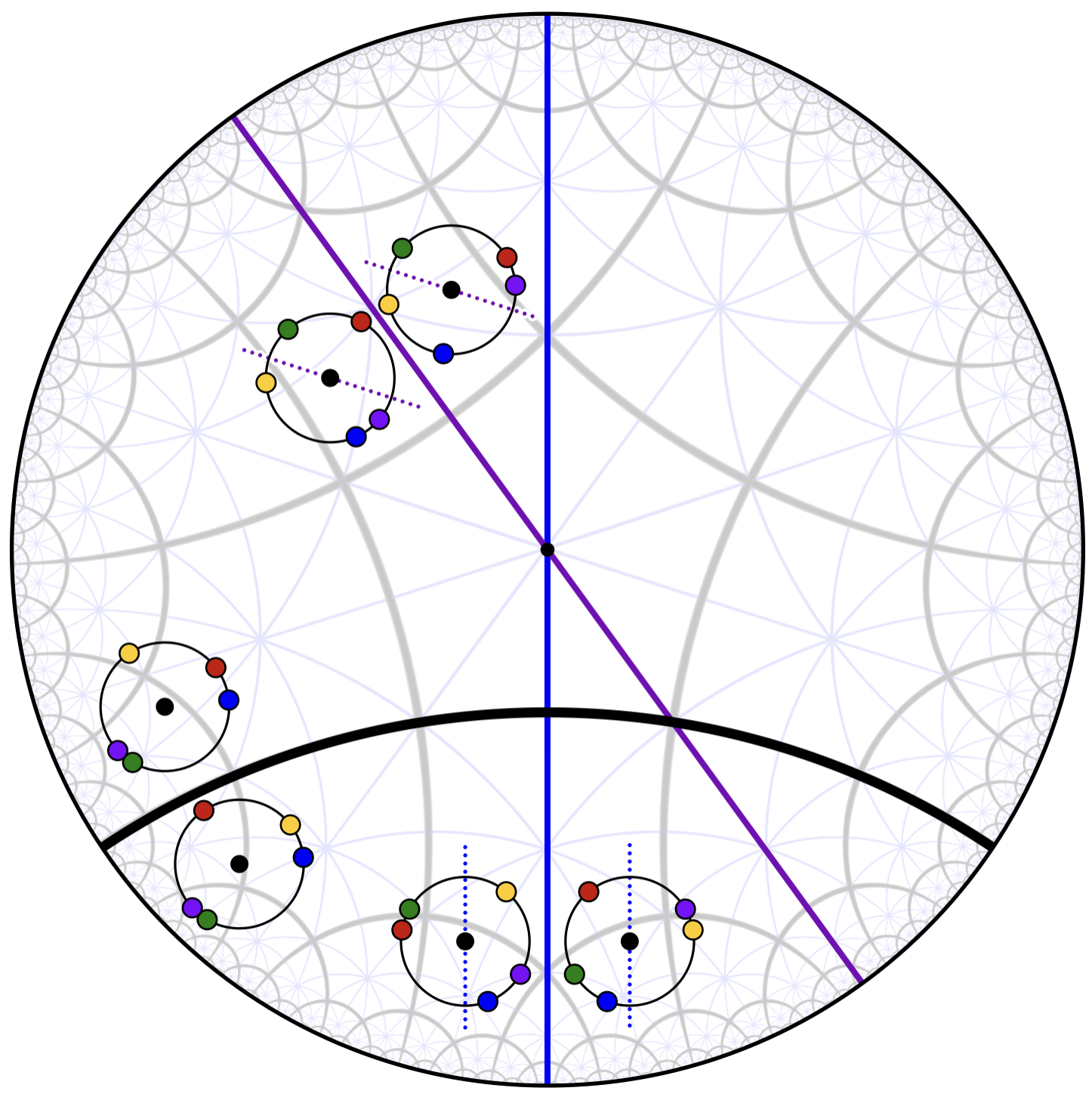}
\begin{picture}(0,0)
\put(-155,302){\footnotesize {$r_0$}}
\put(-155,-5){\footnotesize {$r_0$}}
\put(-244,275){\footnotesize {$r_1$}}
\put(-65,25){\footnotesize {$r_1$}}
\put(-28,65){\footnotesize {$s$}}
\put(-283,65){\footnotesize {$s$}}
\end{picture}

\end{center}
\captionof{figure}{Effect of reflections in certain symmetry axes.}~\label{fig:hypreflections}
\end{figure*}
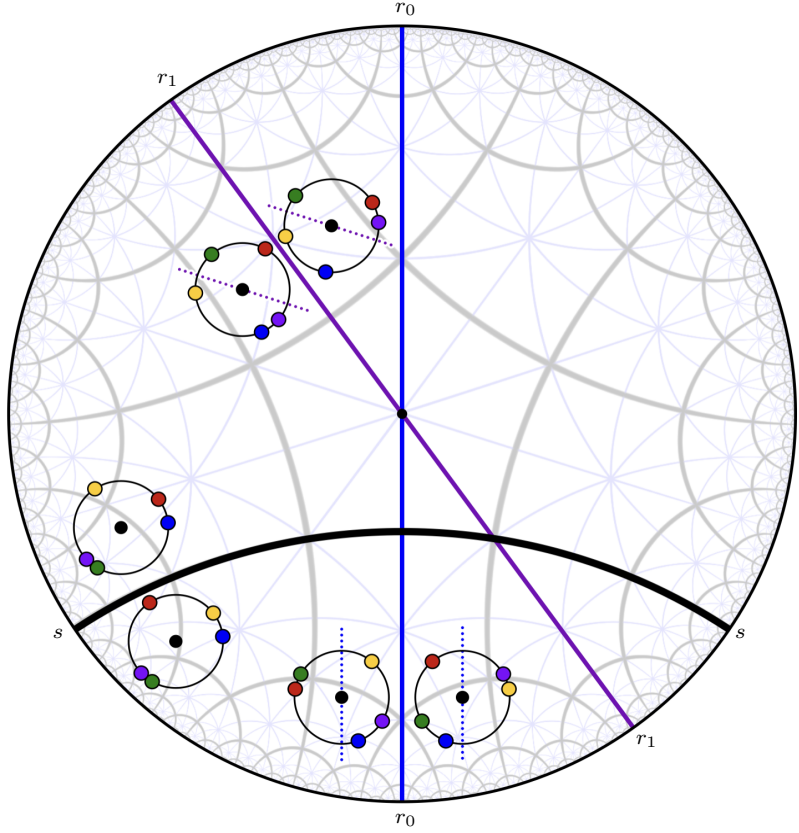

\noindent
In our setup $s$ is a reflection of type~$A$, while the $r_k$ are reflections of type~$B$.
Geometrically, each of them corresponds to a circle inversion (resp. hyperbolic reflection)  with potentially infinite circle radius.
Any other symmetry of the hyperbolic tiling can be expressed by composition of these reflections.
In Figure \ref{fig:hypreflections} we highlighted the three  reflections $s,r_0,r_1$ forming the sides of a triangle that touches the center of the Poincaré disk. 

We now formulate requirements for a group action of $\circ\colon \mathcal{R}\times (S^1)^5\to (S^1)^5$ that are compatibility with our observations about $\mathcal{P}$ and its relation to the tiling.
Each element  in $(S^1)^5$ will be denoted by sequence of letters 
$(a,b,c,d,e)$. The position of these entries corresponds to the sequence {\it blue-purple-red-yellow-green} in our drawings.

\noindent
{{\bf Requirement 1} (Type $A$ symmetries):} 
Crossing the symmetry axis related to $s$ of type $A$ shall result in a flip of 
the yellow and red point. In our drawing for $s$ they have indices~$2$ and~$3$. Thus we set (the bold letters remain fixed):
\[s\circ (\boldsymbol{a},\boldsymbol{b},c,d,\boldsymbol{e}):=(\boldsymbol{a},\boldsymbol{b},d,c,\boldsymbol{e}).\]
\begin{figure}[H]
\centering
\includegraphics[width=0.45\textwidth]{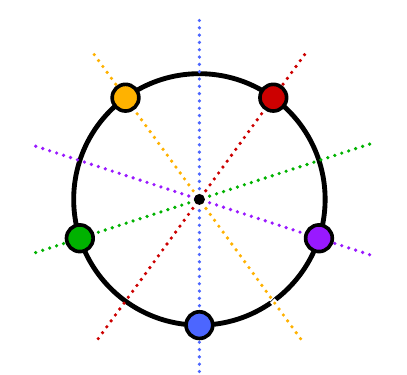}
\begin{picture}(0,0)
\put(-115,203){\footnotesize {$\widetilde{r}_0$}}
\put(-25,132){\footnotesize {$\widetilde{r}_1$}}
\put(-55,25){\footnotesize {$\widetilde{r}_2$}}
\put(-186,25){\footnotesize {$\widetilde{r}_3$}}
\put(-214,132){\footnotesize {$\widetilde{r}_4$}}
\end{picture}
\caption{\label{fig:reflect}The central reflections}
\end{figure}

\noindent
{{\bf Requirement 2} (Type $B$ symmetries):} 
The action for type $B$ symmetries is more complex. 
An axis of type $B$ corresponds to a symmetry that mirrors the 
juzu picture in a certain axes, but at the same time swaps the roles of of two pairs of
points in the juzu, 
Since the juzu image is centered in the origin the reflections correspond to
central axes. 
Figure \ref{fig:reflect} indicates their position.
These reflections are given by
\[
\widetilde{r}_k(z):=-e^{i\cdot 2 k\alpha}\cdot {\overline{z}}
\]
The operation associated to $r_0$ will be a reflection in $\widetilde{r}_0$ followed by swapping the pairs
{\it (red--yellow)} and {\it (green--purple)}. Notice the factor $2$ that doubles the angle of rotation with respect to the $r_k$. The juzu shown in Figure \ref{fig:reflect}  will be a fixed point of each of those operations.
We set
\[
\begin{array}{rcl}
r_0\circ (\boldsymbol{a},b,c,d,e)&:=&\widetilde{r}_0((\boldsymbol{a},e,d,c,b))\\[1mm]
r_1\circ (a,\boldsymbol{b},c,d,e)&:=&\widetilde{r}_1((c,\boldsymbol{b},a,e,d))\\[1mm]
r_2\circ (a,b,\boldsymbol{c},d,e)&:=&\widetilde{r}_2((e,d,\boldsymbol{c},b,a))\\[1mm]
r_3\circ (a,b,c,\boldsymbol{d},e)&:=&\widetilde{r}_3((b,a,e,\boldsymbol{d},c))\\[1mm]
r_4\circ (a,b,c,d,\boldsymbol{b})&:=&\widetilde{r}_4((d,c,b,a,\boldsymbol{b}))\\
\end{array}
\]
Where we consider the application of $\widetilde{r}_k$ to a list element-wise.

\noindent
{{\bf Requirement 3} (Continuity):} 
Furthermore, we should require the map $\Psi$ to be at least continuous. Or even better $C^\infty$ continuous.

\noindent

Figure~\ref{fig:hypreflections} illustrates the actions on the juzus by showing three pairs 
$(\Psi(p),\Psi( g\circ p))$ for the generating reflections $s,r_0,r_1$. The reflection axes for the juzus are indicated in the little circles.
Notice that this group action is compatible with the drawing in Figure~\ref{fig:mapsFull}.

\medskip
We will list a few remarkable consequences 
of our requirements.
We will omit the proofs, since they are quite straight forward, and only require careful bookkeeping of the reflections and index permutations.

\begin{theorem}
The operation $\circ$ defines a group action on $(S^1)^5$.
\end{theorem}
The only thing that has to be shown for this is that $\circ$ respects the defining relations of the symmetry group.
\[
s^2=r_0^2=r_1^2=(sr_0)^2=(sr_1)^4=(r_0r_1)^5=1.
\]
Thus we get a consistent extension of a map $\Psi$ from one triangle of the tiling to the entire hyperbolic plane.

\medskip
It is immediate that our requirements cause the juzus to take certain shapes whenever they lie on a symmetry line. They are compatible with 
the relations of the $(5,4)$-tiling and the pentagon linkage we nailed down so far. 
\begin{figure*}[t]
\noindent
\begin{center}
\includegraphics[width=0.6\textwidth]{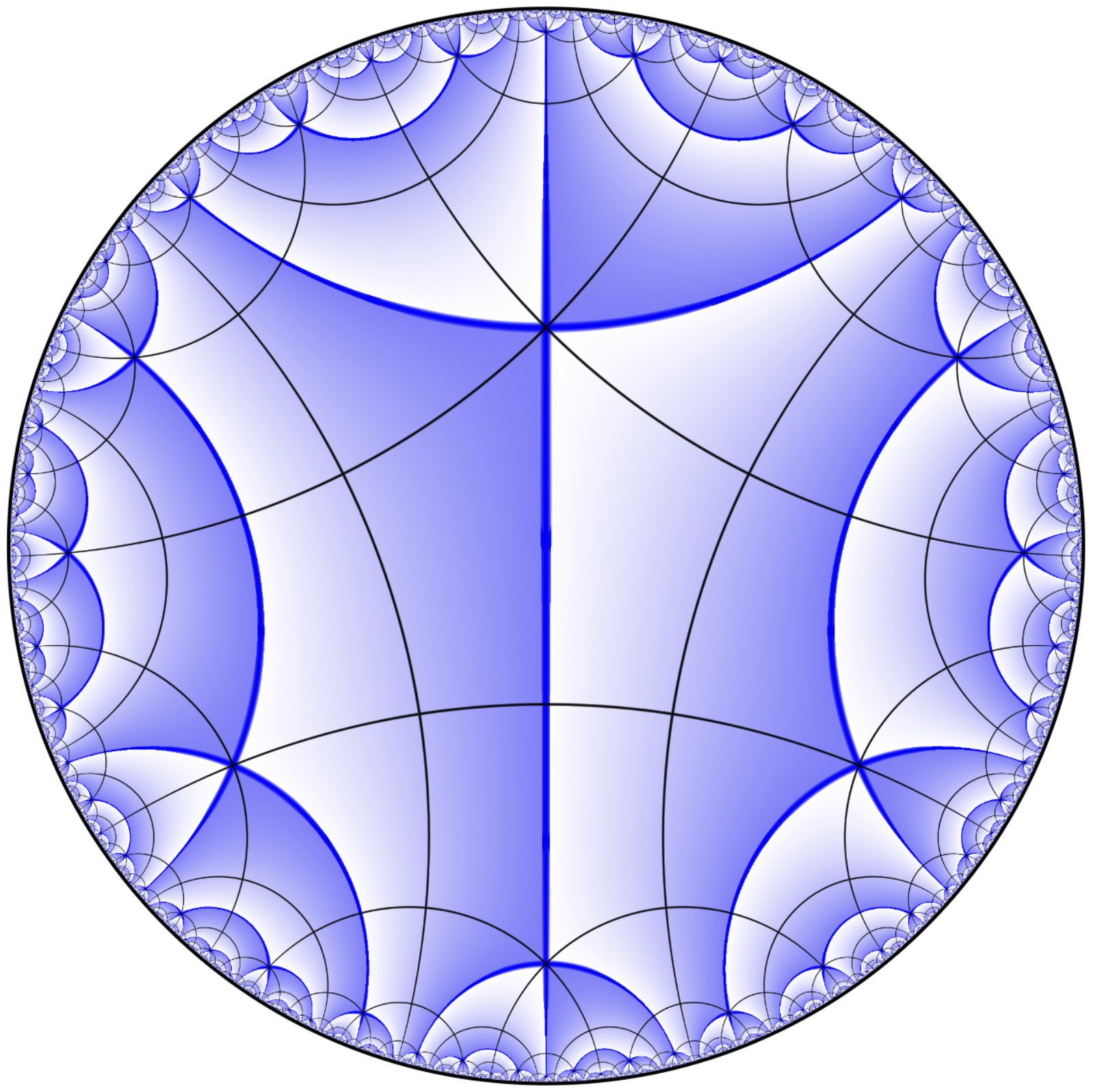}\qquad
\end{center}
\captionof{figure}{A sketch of the behavior of the function $(\Psi(p))_0$ (resp. $\psi$) . Along blue hyperbolic lines this value is constant.}~\label{fig:conf3}
\end{figure*}

\begin{theorem}
Let $\Psi$ satisfy the Requirements 0 to 3. Then
\begin{itemize}
\item[(i)] if $p$ lies on the mirror axis of $s$ then $(\Psi(p))_2=(\Psi(p))_3$,
\item[(ii)] if $p$ lies in the barycenter of one of the pentagons, then 
\[
\begin{array}{rl}
&\big\{(\Psi(p))_k\ \big\vert\ k\in\{0\upto 4\}\big\}\\[2mm]
=&\big\{(\Psi(0))_k\ \big\vert\ k\in\{0\upto 4\}\big\}
\end{array}
\]
\item[(iii)] if $p$ lies on the mirror axis of $r_k$  then $(\Psi(p))_k=-i e^{i2k\alpha}$.
\item[(iv)] Along each mirror axes of type $B$ one of the entries of $\Psi(p)$ is constant.
\end{itemize}
\end{theorem}
%
%

By (iv) each mirror axis of type $B$ may be associated with the color of the index that is unaffected by tis action. 
At each barycenter of a pentagon in the tiling five axes of type $B$ of different color meet. The cyclic order in which these colors meet assumes all different possibilities.
Figure~\ref{fig:conf3} illustrates the position of those axes of type~$B$ for which the value of $(\Psi(p))_0$ remains constant (blue).
These are the blue lines in the drawing. Observe that the blue lines come in different clusters that are disconnected from each other.
If two blue lines cross they must correspond to the same value of  $(\Psi(p))_0$.
The possible values taken at those lines are $(\Psi(p))_0\cdot e^{i\cdot k\alpha}$.
 The color gradient
between the clusters  symbolises the smooth transition from one cluster of lines to the adjacent ones. 
For this transition a rotation of $(\Psi(p))_0$ by $\alpha=72^\circ$ must take place.

\subsection{Just one function}
The next observation reduces our entire problem to finding just one real valued function $\varphi\colon\mathbb{H}\to \mathbb{R}$. Again, the proof is done by careful bookkeeping.

\begin{theorem}
As before let $\Psi$ be a smooth function that satisfies our four requirements. Let $\omega_k=e^{i\cdot k \alpha}$
\[(\Psi(p))_k=
(\Psi(p\cdot \omega_k))_0\cdot(\Psi(0))_k/(\Psi(0))_0
\]
\end{theorem}
%

Let us phrase this proposition in prose:
In order to get the different entries 
of $\Psi(p)$ we have to take the positions at the center
and rotate each of them by the amount
given by $\Psi(p\cdot \omega_k)_0/\Psi(0)_0$. The points $p\cdot \omega_k$
form a centered regular pentagon one of whose vertices is $p$.
Let us define a real valued function $\psi\colon \mathbb{H}\to\mathbb{R}$ 
such that $\Psi(p)_0=e^{i\alpha\psi(p)}$. Then the entire behaviour
of $\Psi$ is reduced to a real valued function $\psi$. This function describes the angle offset of one juzu vector with respect to its position in the center. It takes constant values $\ldots,-3,-2,1,0,1,2,3,\ldots$ for the blue lines in
Figure~\ref{fig:conf3}. Between consecutive stripes the function should interpolate smoothly. Furthermore the function must admit a $180^\circ$ rotation symmetry
around certain edge centres (see below). Every function satisfying these requirements will lead to an admissible parameterisation.  

%
%
%
%
%
\section{If you have a hammer$\ldots$}

$\ldots$everything looks line a nail. What is the best way to get such a function $\psi$?
A natural way to derive such a function is to model it via a conformal (or harmonic) function that 
is consistent with the required symmetries. As mentioned  at the beginning, 
creating visualisations is one of my main research fields and from a different project \cite{MRRG}
we had a very versatile tool to create conformal functions that respect the symmetries of a hyperbolic tiling.
By the Riemann mapping theorem there is a conformal map that maps any topological disk bounded by a Jordan curve
to the unit disk. While it is relatively complicated to calculate conformal mappings for general regions,
there are some fairly good approximation methods for regions bounded by circular arcs. We will use exactly such a mapping
to model the coordinate function $\psi$. 
\begin{figure}[H]
\noindent
\begin{center}
\includegraphics[width=0.13\textwidth]{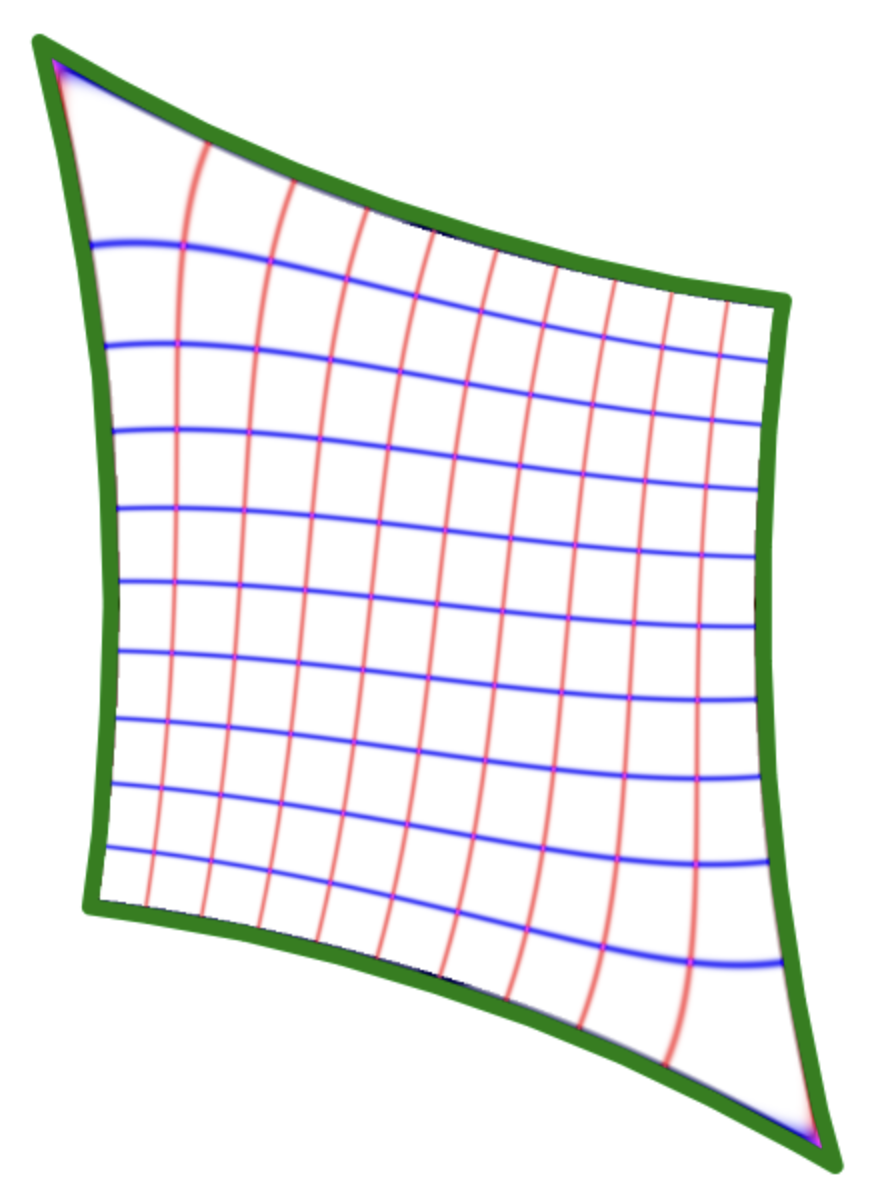}
\includegraphics[width=0.22\textwidth]{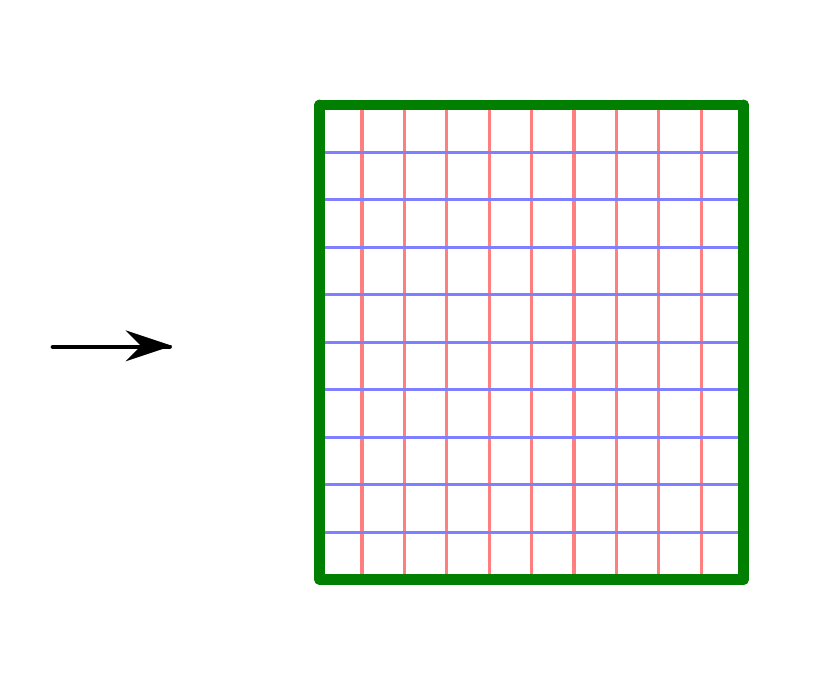}
\end{center}
\captionof{figure}{The conformal mapping from $Q$ to $R$.}~\label{fig:conf5}
\end{figure}

Consider the quadrangle $Q$ bounded by circular arcs meeting in corner angles
$(90^\circ,45^\circ,90^\circ,45^\circ)$. If we want to conformally map this region to a rectangle such that the corners match, then the aspect ratio of the rectangle is not arbitrary. Four points 
on a Jordan curve define a so called {\it conformal modulus}. This is the well defined cross-ratio that arises when the
curve is conformally mapped the unit circle. It turns out that the conformal modulus of $Q$ is
 $m=0.89281029\ldots$.
 This number was approximated by suitable numerical algorithms \cite{Nas22}. A rectangle $R\in \mathbb{C}$  with that aspect ratio can be mapped conformally to $Q$. For later use we assume that the lower left corner of $R$ is the origin and that the sides are parallel to the real and imaginary axes and that its height is 1.
Figure~\ref{fig:conf5} illustrates this map $g\colon Q\to R$. The one additional symmetry requirement that was mentioned above is resembled by the rotation symmetry around the center of $Q$.

\begin{figure}[H]
\noindent
\begin{center}
\includegraphics[width=0.45\textwidth]{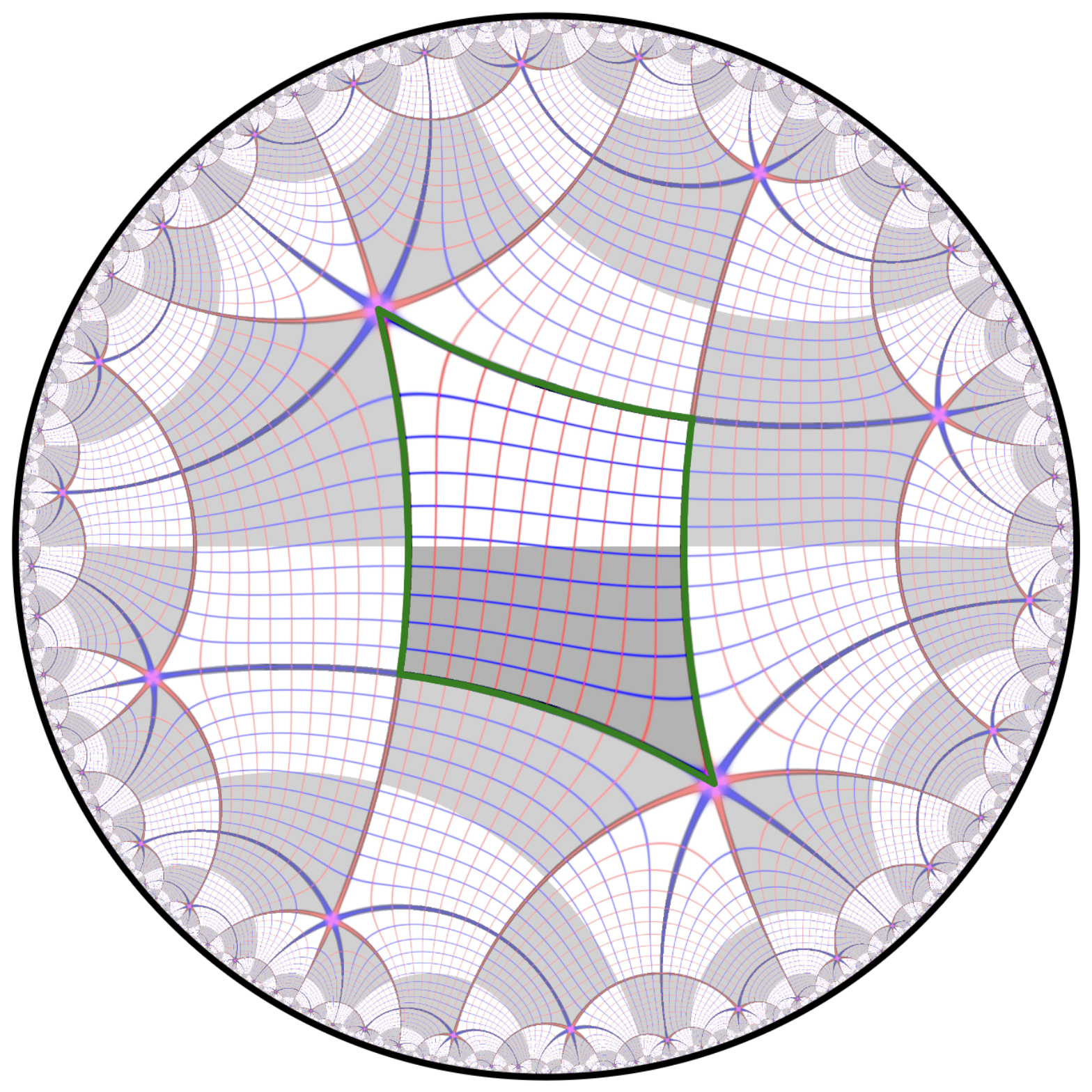}\qquad
\end{center}
\captionof{figure}{The pentagon linkages at the centres of the cells in each dimension.}~\label{fig:mapping}
\end{figure}

Why did we take this particular region $Q$? For an answer consider Figure~\ref{fig:mapping}.
In that image we overlayed two symmetry groups. In the background
you see a black and white checkerboard pattern of the $(5,4)$-tiling. To emphasise the symmetry of the region $Q$
we have applied a hyperbolic transformation to the tiling that aligns one side of a pentagon with the real axis.
On top of that one can see (with a green border) the tile $Q$. Notice that it has exactly the same hyperbolic area as the pentagons in the tiling (one can dissect the pentagons along their symmetry axis and reassemble them to get $Q$).
Now one can take the boundaries of $Q$ (these are hyperbolic straight lines) and consider the reflection
group generated by taking them as mirrors. Since the corner angles of $Q$ are divisors of $180^\circ$ they generate a 
reflection group in the hyperbolic plane. In other words: like the right-angled pentagons the 
tile $Q$ seamlessly covers the hyperbolic plane.

Now a powerful tool  of conformal mapping theory comes into play: The {\it Schwarz Reflection Principle}.
It states that if a conformal map $f \colon D\to W$ between two topological discs is given such that 
the region $D$ has a circular arc $d$ as part of its boundary and this arc is mapped to a circular arc $w$ bounding part of $W$,
then the map $f$ can be extended in the following way. Let $r_d$ and $r_w$ denote the circular reflections
in the circles supporting $d$ and $q$, respectively. Let $p\not\in D$ be a point with $r_d(p)\in D$ then its image
can be defined by $p\mapsto r_w(f(r_d(p)))$. This map conformally extends $f$ across the boundaries of $D$ and $W$.

Our tiling created by $Q$ together with the conformal map $g\colon Q\to R$ can be extended consistently to become a map
$g\colon \mathbb{H}\to\mathbb{C}$. The fact that at each corner of the tiling by $Q$ either $4$ or $8$ tiles come together
ensures that this map is in fact surjective and single valued. Its inverse is multivalued and exhibits monodromy. 
Observe that the function $g$ has a constant imaginary-value along the blue lines shown in Figure~\ref{fig:mapping}.
In particular this imaginary-value is constant at certain symmetry axes of the pentagons of the $(5,4)$ tiling. For each pentagon there is exactly one symmetry axis where the imaginary part of $g$ is constant. Since we assumed the rectangle $R$ to be aligned
with the coordinate axis and  it has height 1 the function $g$ takes exactly integer values along these symmetry axes where it is constant. The imaginary part $\mathbf{im}(g(p))$ has exactly the desired properties, to behave like the function $\psi(p)$ we are looking for.

\begin{figure}[H]
\noindent
\begin{center}
\includegraphics[width=0.48\textwidth]{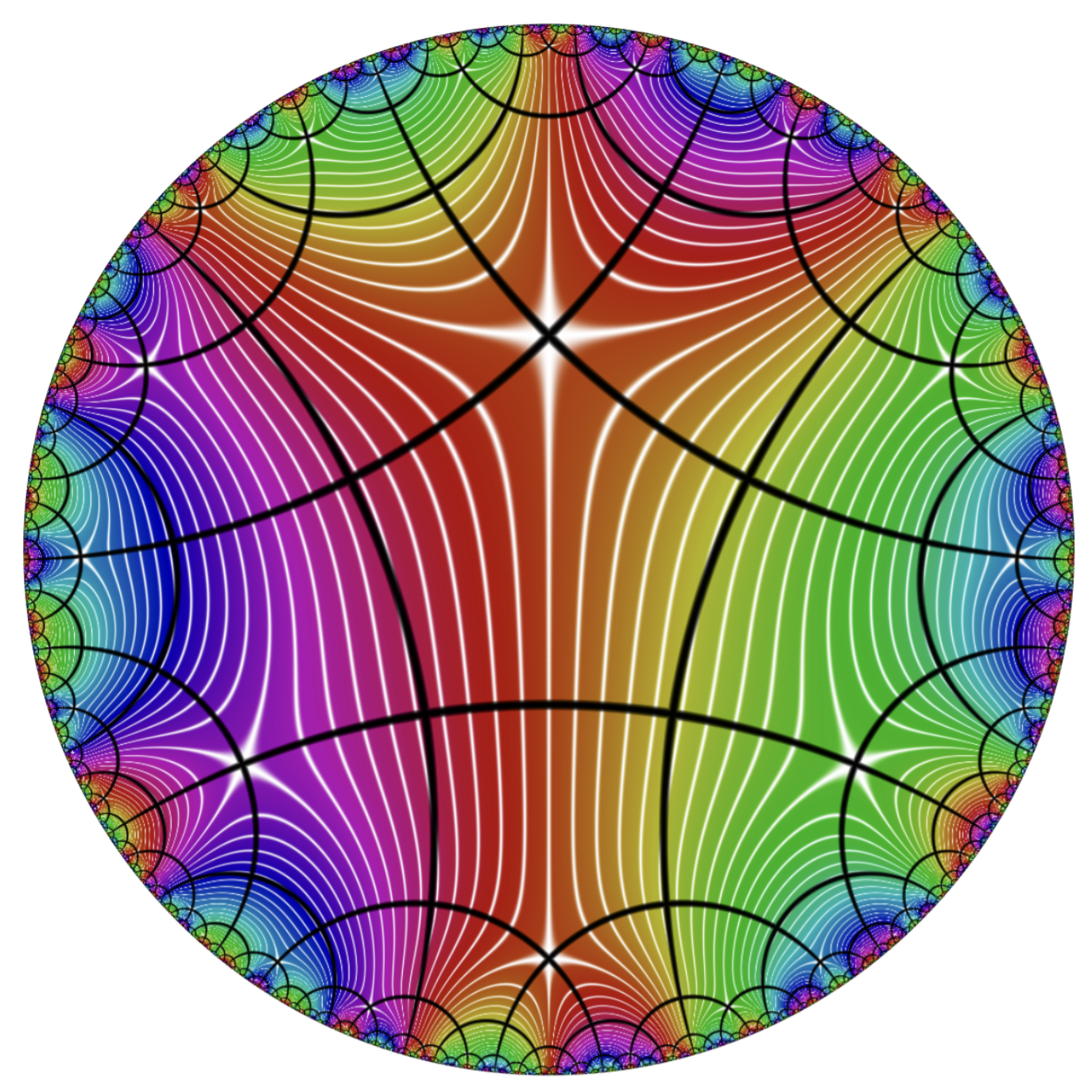}\qquad
\end{center}
\captionof{figure}{The pentagon linkages at the centres of the cells in each dimension.}~\label{fig:flow}
\end{figure}

Figure \ref{fig:flow} 
illustrates the {\it flow} of the $y$-coordinate  in the hyperbolic disk. This time the
hyperbolic tiling is again moved to the position we used before where one of the coordinate pentagons is centred at the origin.
Values in the interval $[0,5]$ are indicated by a rainbow hue with red resembling $0$. The white lines are lines at which the 
function $g$ takes a constant value.
Comparing this image to Figure~\ref{fig:conf3} exemplifies that we can take $\psi(q):=\mathbf{im}(g(p))$
in our coordinate function 
\[\Psi(p)_0=e^{i\alpha \psi(p)}.\] 
This is exactly a function we were looking for.
Clearly, it is not the only such interpolating function but it is one that is particularly simple to construct, provided that
the function~$g$ that maps the hyperbolic quadrilateral $Q$ to $R$ is given.

\section{The cooking recipe}
Let us stitch everything together to form a kind of cooking recipe (one might call it algorithm, if there were not many
nasty implementation details) of how to get from a point $q\in \mathbb{H}$ to a concrete
juzu or pentagon linkage.

\begin{figure}[H]
\noindent
\begin{center}
\includegraphics[width=0.4\textwidth]{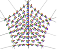}\\
\end{center}
\captionof{figure}{The global picture of the parameterisation.}~\label{fig:film}
\end{figure}

\begin{itemize}
\item[1:] Take a point $q_0\in \mathbb{H}$ and calculate the points $q_k=q_0\cdot e^{i\cdot k \alpha}$.
\item[2:] Calculate the values $\psi_k = \mathbf{im}(\psi(q_k))$.
\item[3:] Calculate the vectors $v_k = e^{i\cdot \psi_k\alpha}\cdot e^{i\cdot 2k\alpha}$.
\item[4:] To the  $v_k$ apply Springborn's normalisation to get $v'_k$ with $\sum v'_k=0$.
\item[5:] Set $p_0=0$ and $p_k=p_{k-1}+v'_k$.
\end{itemize}
The resulting points $p_i$ form a pentagonal linkage smoothly controlled by $q$ and 
respecting all our symmetry requirements.
Figure~\ref{fig:film} shows the result of this procedure for the central pentagon in the $(5,4)$-tiling. A continuous implementation can be found at the webpage 
 \href{https://science-to-touch.com/Visu/Pentagon/Pentagon_Minimal.html}{\tt math-visuals.org/Pentagon}.
 
\begin{figure*}[t]
\noindent
\begin{center}
\includegraphics[width=0.9\textwidth]{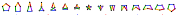}\\[-1mm]
\includegraphics[width=0.9\textwidth]{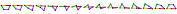}\\[-1mm]
\includegraphics[width=0.9\textwidth]{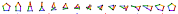}\\[-1mm]
\end{center}
\captionof{figure}{The global picture of the parameterisation.}~\label{fig:film2}
\end{figure*}

The rotations of the pentagrams after the Springborn normalisation  also have been chosen in a way that the satisfy our Requirements 0 to 3 on the level of pentagram linkages. The rows shown in Figure~\ref{fig:film2} are sequences of stills of particular paths of the control point in $\mathbb{H}$.

\section{Some conjectures based on numerical data}
Let us end with a few remarkable properties that seem to be satisfied by our 
function $g\colon \mathbb{H}\to \mathbb{C}$. 
In our cooking recipe we evaluate $g$ at the regular pentagon points  
$q_k=q\cdot e^{i\cdot k \alpha}$. There is strong numerical evidence
that the sum
$\sum_k g(q_k)$ is constant and takes the value ${5\over  2}m$. In particular, its imaginary part (the part that we use) is zero. 
This translates to a fact with  physical relevance. The angle sum of the juzu offsets
remains zero. This property can be carried over 
through the Springborn normalisation.
Thus if we consider a physical realisation of the linkage and give it an initial push that triggers a movement with zero angular momentum, its concrete trajectory can
be parameterised by a path of our control point. 
The lower row of Figure~\ref{fig:film} shows a remarkable movement that demonstrates that while a linkage can remain angular momentum zero it can perform a rotation by $72^\circ$ in \cite{TD12} such movements are studied. There also a remarkable connection is made to the fact that a cat can land on its feet when falling from the second floor.

Another remarkable property that can be numerically observed is relating to conformality of our map.
If we consider
\[
q\mapsto (\mathbf{im}(g(q_1)),\ldots,\mathbf{im}(g(q_5))),
\]
it maps the hyperbolic plane $\mathbb{H}$ to some 2-dimensional
surface in $\mathbb{R}^5$. This map appears to be conformal.



\medskip
\bibliographystyle{plain} 
{\small

}
\medskip

\noindent

\noindent
{\footnotesize
\noindent
{\sc Department of Mathematics, Technical University of Munich, Germany}\\
Email address: {\tt richter@tum.de}
}

\end{multicols}
\end{document}